\documentclass[12pt,a4paper]{article}

\usepackage{amssymb}

\def\cM{{\cal M}}
\def\cT{{\cal T}}
\def\oD{{\overline D}}
\def\oA{{\bar {\bf A}}}
\def\oM{{\overline{\cal M}}}
\def\oC{{\overline{\cal C}}}
\def\qed{{\hfill $\diamondsuit$}}
\def\tk{{\tilde k}}
\def\tDIV{{\widetilde {\rm DIV}}}
\def\DIV{{\rm DIV}}
\def\cL{{\cal L}}
\def\CP{{{\mathbb C}{\rm P}}}

\def\Z{{\mathbb Z}}
\def\ident{
\begin{picture}(18,10)
\put(3,-1){$\rightarrow$}
\put(3.5,2.5){$\sim$}
\end{picture}
}

\usepackage{theorem}

\newtheorem{theorem}{Theorem}
\newtheorem{proposition}{Proposition}[section]

\newtheorem{lemma}[proposition]{Lemma}

{\theorembodyfont{\rmfamily}
\newtheorem{definition}[proposition]{Definition}

\newtheorem{remark}[proposition]{Remark}

}

\include{epsf}

\title{Intersection numbers with Witten's top Chern class}

\author{Sergei Shadrin\thanks{Department of Mathematics, 
Stockholm University,
SE-106 91 Stockholm, Sweden and  Moscow Center for Continuous
Mathematical Education, B.~Vlassievskii per.~11, 119002 Moscow,
Russia. E-mails: shadrin@math.su.se, shadrin@mccme.ru.
Partly supported by the grants RFBR-05-01-01012a, NSh-1972.2003.1,
NWO-RFBR-047.011.2004.026 (RFBR-05-02-89000-NWO-a), 
by the G{\"o}ran Gustafsson foundation, and by
Pierre Deligne's fund based on his 2004 Balzan prize in mathematics.},
Dimitri Zvonkine\thanks{
Institut math{\'e}matique de Jussieu,
Universit{\'e} Paris~VI, 175, rue du Chevaleret,
75013 Paris, France. E-mail: zvonkine@math.jussieu.fr.
Partly supported by the ANR project ``Geometry and
Integrability in Mathematical Physics''.}}

\date{\today}

\begin{document}

\maketitle

\begin{abstract}
Witten's top Chern class is a particular cohomology class on the
moduli space of Riemann surfaces endowed with $r$-spin structures. 
It plays a key role in Witten's conjecture relating to the intersection
theory on these moduli spaces.

Our first goal is to compute the integral of Witten's class over
the so-called double ramification cycles in genus~1. We obtain
a simple closed formula for these integrals.

This allows us, using the methods of~\cite{Shadrin}, to find
an algorithm for computing the intersection numbers of the
Witten class with powers of the $\psi$-classes
(or tautological classes) over any moduli space of $r$-spin
structures, in short, all numbers involved in Witten's conjecture.
\end{abstract}

\section{Introduction}
\label{Sec:Intro}

\subsection{Aims and purposes}
\label{Ssec:Aims}

In 1991 E.~Witten formulated two conjectures relating to the 
intersection theory of moduli spaces of 
curves~\cite{Witten1, Witten2}, motivated by two dimensional
gravity.

The first conjecture involves moduli spaces of
stable curves and the tautological 2-cohomology classes
on them (also called $\psi$-classes). The intersection
numbers of powers of the $\psi$-classes can be arranged into a
generating series that is claimed to be a solution of
the Korteweg -- de~Vries (or KdV) hierarchy of partial
differential equations. This conjecture was first proved by
M.~Kon\-tse\-vich in~\cite{Kontsevich}. At present there are
several alternative 
proofs:~\cite{OkoPan},~\cite{Mirzakhani},~\cite{KazLan},~\cite{KimLiu}.

The second conjecture is still open, and even giving a precise 
formulation required joint work by several people (references
are given below).
It involves a more complicated moduli space, called the
{\em space of\/ $r$-spin structures}. Apart from the $\psi$-classes,
one considers one more cohomology class, called the {\em Witten
top Chern class}, or just {\em Witten's class} for shortness. 
We are interested in the intersection numbers of
Witten's class with powers of the $\psi$-classes.
These intersection numbers can, once again, be arranged into
a generating series, and this series is claimed to give a
solution of the $r$th Gelfand-Dikii (or the $r$-KdV) hierarchy.

The precise definition of Witten's top Chern class is rather
involved; however it is known to satisfy quite simple
factorization rules. 

In~\cite{Shadrin} the first author found an ``almost-algorithm''
for computing some of the intersection numbers arising in Witten's second
conjecture {\em using only the factorization rules for the
Witten class}. More precisely, the factorization rules allow
one to express more complicated intersection numbers
via simpler ones, until one arrives at unsimplifiable 
cases. These can be of two types. 
(i)~Integrals of Witten's class (with no $\psi$-classes)
over genus~$0$ moduli spaces. These numbers are well-known.
(ii)~Integrals of Witten's class (with no $\psi$-classes)
over some special divisors on genus~$1$ moduli space.
(These divisors have a rather cumbersome name of {\em double
ramification divisors} - see below.)
When numbers of second type appeared in the course of computations
the algorithm blocked without giving an answer.

The purpose of this note it twofold. 

First, compute the integrals of Witten's class over the double 
ramification divisors in genus~1. 
It turns out that a simple closed formula exists
for these integrals. 

Second, complete and give a coherent exposition of
the algorithm for computing Witten's intersection
numbers. Our computation uses only factorization
rules for Witten's class. Therefore, we obtain the following theorem.

\begin{theorem} \label{Thm:algorithm}
The intersection numbers of Witten's class with powers of the
$\psi$-classes are entirely determined by (i) genus~$0$ intersection
numbers involving no $\psi$-classes, and (ii) the factorization 
rules for Witten's class.
\end{theorem}

\subsection{Main definitions}
\label{Ssec:Definitions}

\subsubsection{Moduli spaces.}
$\cM_{g,n}$ is the moduli space of smooth complex genus~$g$
curves with $n \geq 1$ distinct numbered marked points.
$\oM_{g,n}$ is its Deligne-Mumford compactification,
in other words, the moduli space of stable curves. Over $\oM_{g,n}$ we
define $n$ holomorphic line bundles $\cL_i$. The fiber of
$\cL_i$ over a point $a \in \oM_{g,n}$ is the cotangent
line to the corresponding stable curve $C_a$ at the
$i$th marked point. The first Chern classes $\psi_i = c_1(\cL_i)$
of these line bundles are called the {\em tautological classes}.

\subsubsection{Spaces of $r$-spin structures.}
Chose an integer $r \geq 2$ and pick $n$ integers
$a_1, \dots a_n \in \{ 0, \dots, r-1 \}$ in such a way
that $2g-2 - \sum a_i$ is divisible by~$r$. The numbers
$a_1, \dots, a_n$ are assigned to the marked points
$x_1, \dots, x_n$. On a smooth
curve~$C$ one can find $r^{2g}$ different line bundles $\cT$ with an
identification
$$
\cT^{\otimes r} \ident K \left(-\sum a_i x_i\right).
$$
The space of smooth curves endowed with such a line bundle
$\cT$ is called the {\em space of $r$-spin structures}
and denoted by $\cM^{1/r}_{g; a_1, \dots, a_n}$. It is an unramified
$r^{2g}$-sheeted covering of $\cM_{g,n}$.

A compactification of this space, denoted by 
$\oM^{1/r}_{g;a_1, \dots, a_n}$, was constructed  
in~\cite{Jarvis} and~\cite{AbrJar} (see also~\cite{Chiodo2}
for a slightly simplified version).
It is a smooth orbifold (or stack), and there is a finite
projection mapping
$$
p:\oM^{1/r}_{g;a_1, \dots, a_n} \rightarrow \oM_{g,n}.
$$
The construction uses the so-called {\em Jarvis-Vistoli twisted
curves}, i.e., curves that are themselves endowed with an
orbifold structure. The stabilizers of the marked points and the nodes
are equal to $\Z/r\Z$, while the stabilizer of any other point is
trivial. $\cT$ is then an $r$th root of $K(-\sum a_i x_i)$
{\em in the orbifold sense}. Alternatively, we can forget about the
orbifold structure of the curve and consider only the sheaf of 
invariant sections of~$\cT$. We then obtain a rank one 
torsion-free sheaf rather than a line bundle.

\subsubsection{Witten's class.}
The rank one torsion-free sheaf (of invariant sections of)
$\cT$ is defined on the 
universal curve $\oC^{1/r}_{g;a_1, \dots, a_n}$ over
$\oM^{1/r}_{g;a_1, \dots, a_n}$. Consider its push-forward
to the space $\oM_{g;a_1, \dots, a_n}$ itself. First assume
that for each curve $C$ we have $H^0(C, \cT) = 0$. 
Then the spaces $H^1(C,\cT)$ form a vector bundle $V^\vee$
over $\oM^{1/r}_{g;a_1, \dots, a_n}$. We denote by $V$
the dual vector bundle and define Witten's class as
$$
c_W (a_1, \dots, a_n)  = c_W = \frac{1}{r^g} \; p_* \, c_{\rm top} (V).
$$
In other words: take the top Chern class of $V$,
push it from $\oM^{1/r}_{g; a_1, \dots, a_n}$ to $\oM_{g,n}$,
and divide by $r^g$. By the Riemann-Roch formula, the (complex)
degree of Witten's class is
$$
\deg c_W = \frac{(r-2)(g-1) + \sum a_i}r .
$$

Unfortunately, in general $\cT$ has both
$0$- and $1$-cohomologies. The definition of Witten's
class $c_W$ in this case is much more involved. There exist
two algebro-geometric constructions: see~\cite{PolVai} 
and~\cite{Chiodo1}. (Note that the definition of
$c_W$ is special to our situation and uses the identification
of $\cT^{\otimes r}$ with the canonical bundle. No
general constructions from algebraic geometry are expected
to work.)

Witten's class satisfies the following vanishing property, that
we will use as an axiom: 

{\bf If one of the $a_i$'s equals $r-1$ then $c_W = 0$}.

\subsubsection{Factorization rules.}
\begin{figure}[h]
\begin{center}
\
\epsfbox{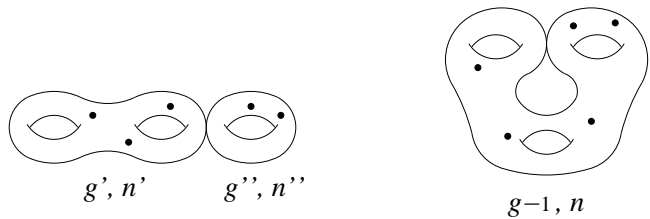}
\caption{Two possible degenerations of a stable curve.}
\label{Fig:degeneration}
\end{center}
\end{figure}

We are interested in the restriction of Witten's class to
the boundary components of the moduli space $\oM_{g,n}$.
There are two types of boundary components 
(see Figure~\ref{Fig:degeneration}): those isomorphic
to $\oM_{g',n'+1} \times \oM_{g'',n''+1}$, $n'+n''= n$, $g'+g''=g$
and the unique component isomorphic to 
$\oM_{g-1, n+2}/{\mathbb Z}_2$.

In the first case, assume for simplicity that the 
the marked points $x_1, \dots, x_{n'}$ are on the first
component of the curve, while $x_{n'+1}, \dots, x_n$ are 
on the second component. There is a unique choice of 
$a', a'' \in \{ 0, \dots, r-1 \}$ such that
$$
2g'-2 - a' - \sum_{i=1}^{n'} a_i 
\qquad  \mbox{and}  \qquad 
2g''-2 - a'' - \sum_{i=n'+1}^n a_i
$$
are both divisible by~$r$. We have $a'+a'' = r-2$ or 
$a'=a''=r-1$. For the second type of boundary component, we have to sum over
all choices of $a', a''$ such that $a'+a'' = r-2$.
Now we can formulate the factorization rules. 

{\bf The restriction of Witten's class to
the first type boundary component equals 
$$
c_W(a_1, \dots, a_n) = 
c_W(a_1, \dots, a_{n'}, a') \times c_W(a_{n'+1}, \dots, a_{n}, a'').
$$

The restriction of Witten's class to
the second type boundary component equals 
$$
c_W(a_1, \dots, a_n) = \frac12 
\sum_{a'+a'' = r-2}
c_W(a_1, \dots, a_n, a', a'').
$$
}

The vanishing property and the factorization rules (for the
Polishchuk-Vaintrob construction) are proved in~\cite{Polishchuk}.

\subsubsection{Intersection numbers.}

We use the standard notation for the intersection numbers of Witten's
class with powers of the $\psi$-classes:
$$
\left< \tau_{d_1, a_1} \dots \tau_{d_n, a_n} \right>
= \int\limits_{\oM_{g,n}}
c_W(a_1, \dots, a_n) \, \psi_1^{d_1} \dots \psi_n^{d_n}.
$$
Although the genus $g$ is determined by the elements of the bracket
by
$$
(r+1)(2g-2+n) = \sum (rd_i + a_i +1),
$$
we will sometimes recall it in a subscript. The integer $r \geq 2$
is supposed to be fixed once and for all throughout the paper.

Since double ramification cycles on genus~$1$ moduli spaces
will play a special role, let us introduce them here (this is
a particular case of Definition~\ref{Def:DR-cycle}). Choose
$n$ integers $k_1, \dots, k_n$ satisfying $\sum k_i = 0$. We
assume that at least one of the $k_i$'s is different from~$0$.

We will usually assume that the list $( k_1, \dots, k_n )$
starts with the positive integers and ends with the negative ones,
the zeroes being in the middle. We will sometimes use the notation
$$
( k_1, \dots, k_{n_+} \; | \; 0, \dots, 0 \; | \; \tk_1, \dots, \tk_{n_-})
$$
with only nonnegative integers instead of
$$
(k_1, \dots, k_n) = 
(k_1, \dots, k_{n_+}, 0, \dots, 0,  -\tk_1, \dots, -\tk_{n_-}).
$$

To this list of integers assign the set $D(k_1, \dots, k_n)$
of smooth genus~1 curves $(C, x_1, \dots, x_n)$ 
such that $\sum k_i x_i$ is the divisor of a function.
Let the cycle $\oD(k_1, \dots, k_n)$ be the closure of 
$D(k_1, \dots, k_n)$.

\begin{definition} \label{Def:divisor}
We call  $\oD(k_1, \dots, k_n)$ a {\em double ramification cycle}.
The integral of Witten's class over this cycle is denoted by
$$
\int\limits_{\oD_{k_1, \dots, k_n}} \!\!\!\!\!\!\! c_W (a_1, \dots, a_n)
=
\left<
\begin{array}{ccc}
k_1 & \dots & k_n \\
a_1 & \dots & a_n
\end{array}
\right>
=
$$
$$
=
\left<
\begin{array}{ccc|ccc|ccc}
k_1 & \dots & k_{n_+} & 0 & \dots & 0 & \tk_1 & \dots & \tk_{n_-} \\
a_1 & \dots & \dots & \dots & \dots & \dots & \dots & \dots & a_n
\end{array}
\right>
$$
\end{definition}

\begin{theorem} \label{Thm:formula}
We have
$$
\left<
\begin{array}{ccc}
k_1 & \dots & k_n \\
a_1 & \dots & a_n
\end{array}
\right>
= \left( \frac12 \sum_{i=1}^n k_i^2 - 1 \right) \cdot \frac1{24}
 \frac{(n-1)!}{r^{n-1}} \prod_{i=1}^n (r-1-a_i).
$$
\end{theorem}

\subsection{Acknowledgements}
This work was completed during the second author's visit
to the Stockholm University. We would like to thank the
University for its hospitality and the members of the
mathematical seminar for their interest.

\section{Preliminaries} 
\label{Sec:Prelim}

\subsection{Admissible coverings}
\label{Ssec:Adm}

Consider a map $\varphi$ from a smooth curve $C$ to the sphere $S =
\CP^1$. Mark
all ramification points on $S$ and all their preimages on~$C$.
Now choose several disjoint simple loops on $S$, that do
not pass through the marked points. Suppose that if we contract
these loops we obtain a stable genus~0 curve $S'$. Now contract
also all the preimages of the loops in~$C$ to obtain a
nodal curve $C'$ that turns out to be automatically stable. 
We have obtained a map $\varphi'$ from a nodal curve 
of genus~$g$ to a stable curve of genus~0. It has the
same degree over every component of $S'$. Moreover, at each
node of $C'$, the projection $\varphi'$ has the same local
multiplicity on both components meeting at the node.

\begin{definition} \label{Def:AdmCov}
A map from a stable curve of genus~$g$ to a stable curve of
genus~0 topologically equivalent to a map described above is called
an {\em admissible covering}. 
\end{definition}

The space of all admissible coverings
with prescribed ramification types over the marked points is
very useful for the study of moduli spaces (see~\cite{Ionel}). 
It is not normal but can be normalized, 
and thus one can study its intersection theory.
We refer to Ionel's work~\cite{Ionel} for detailed definitions.
This construction can be slightly extended by allowing additional marked
points on the curve~$C$. Their images on the genus~0 curve will
also be marked (although they are not ramification points). All
the other preimages on~$C$ of these new marked points on~$S$
should also be marked.

We will be particularly interested in the space of admissible
covering with multiple ramifications only over~2 points labeled
with $0$ and~$\infty$, the other ramification points being simple.

\begin{definition} \label{Def:DR-Space}
Consider the space of admissible coverings of some given genus~$g$
with prescribed ramification types $(k_1, \dots, k_{n_+})$
and $(\tk_1, \dots, \tk_{n_-})$, $\sum k_i = \sum \tk_i$, 
over two points labeled $0$ and $\infty$, and with
simple ramifications elsewhere. We also suppose that
there are $n_0$ additional marked points on $C$.
The normalization of this space
is called a {\em double ramification space} or a {\em DR-space}. 
It is denoted by 
$$
\oA = \oA(k_1, \dots, k_{n_+}, 
\overbrace{0, \dots, 0}^{n_0}, -\tk_1, \dots, -\tk_{n_-}).
$$
We do not include the genus~$g$ in our notation, although,
of course, the space of coverings does depend on~$g$.
\end{definition}

If $N$ is the total number of marked points on the curve~$C$, we can
consider the forgetful map $j: \oA \rightarrow \oM_{g,N}$ that
forgets the admissible covering retaining only its source curve.
Since the curve~$C$ is automatically stable, the map~$j$ is actually
an injection and an isomorphism with its image. The pull-backs by
$j$ of the $\psi$-classes on $\oM_{g,N}$ coinside with the
$\psi$-classes naturally defined on~$\oA$.

Another forgetful map $f: \oA \rightarrow \oM_{0,n+2g-2}$,
$n = n_++n_-+n_0$ takes an admissible to its image genus~0
curve. This map satisfies the following crucial property.

\begin{lemma} [Ionel's lemma, see~\cite{Ionel}] \label{Lem:Ionel}
The map $f$ sends the fundamental homology class of $\oA$
to (a multiple of) the fundamental homology class of
$\oM_{0,n+2g-2}$. We have
$$
\psi_i(\oA) = \frac1{k_i} \psi_0(\oM_{0,n+2g-2}),
\qquad
\psi_{n_++n_0+i}(\oA) = \frac1{\tk_i} \psi_\infty(\oM_{0,n+2g-2}).
$$
In other words, a $\psi$-class on $\oA$ at a marked zero
or pole of the admissible covering coincides, up to a constant, 
with the pull-back of the $\psi$-class at $0$ or $\infty$
of the genus~0 moduli space.
\end{lemma}

Finally, let us define double ramification cycles.
For a given moduli space $\oM_{g,n}$ choose an
integer $p$, $0 \leq p \leq g$  and $n+p$ integers $k_i$. We suppose that
$\sum k_i = 0$ and that none of the $k_{n+1}, \dots, k_{n+p}$
vanishes.

\begin{definition} \label{Def:DR-cycle}
Consider the set of smooth curves $(C, x_1, \dots, x_n) \in \cM_{g,n}$
such that there exist $p$ more marked points $x_{n+1}, \dots, x_{n+p}$
and a meromorphic function on $C$ with no zeroes or poles outside of
$x_1, \dots, x_{n+p}$, the orders of zeroes or poles 
being prescribed by the list
$k_1, \dots, k_{n+p}$ ($k_i >0$ for the zeroes, $k_i <0$ for the
poles, and $k_i=0$ for the marked points that are neither zeroes nor
poles). The closure of this set in $\oM_{g,n}$
is called the {\em double ramification
cycle} or a {\em DR-cycle}.
\end{definition}

By a generalization of Mumford's argument in~\cite{Mumford},
one can show that the 
codimension of a DR-cycle is equal to $g-p$ whenever there is at
least one positive and one negative number among $k_1, \dots, k_n$. 
Assuming that this condition is satisfied we see that
for $p=g$ the DR-cycle coincides with the moduli space $\oM_{g,n}$.
Faber and Pandharipande~\cite{FabPan} proved
that the cohomology classes Poincar{\'e} dual to
any DR-cycle belongs to the tautological ring of the moduli space of 
curves. Their proof may be used to obtain new relations for 
intersection numbers with Witten's class.

It follows that 
the forgetful map $h:\oA(k_1, \dots, k_{n+p}) \rightarrow \oM_{g,n}$
sends the fundamental comology class of $\oA(k_1, \dots, k_{n+p})$
to (a multiple of) the fundamental homology class of the
DR-cycle.

\paragraph{Conventions.}
Figure~\ref{Fig:diagramma} represents a DR-space
$$
\oA(k_1+1, \dots, k_{n_+}, 0, \dots, 0, -\tk_1, \dots,
-\tk_{n_-},-1).
$$
It also shows the corresponding maps $j$, 
$f$ and $h$. In this figure, as well as in the subsequent
figures and in the text we follow the following conventions.
A {\bf cross} represents a critical point or a ramification point 
of an admissible covering. {\bf Round black dots} represent
the marked points $x_1, \dots, x_n$ (they are not forgotten
under the map~$h$). The images of these points under the maps
$j$, $f$, and $h$ are also represented as round black dots.
{\bf Square black dots} represent the points $x_{n+1}, \dots, x_{n+p}$
(they are forgotten under the map~$h$). Finally,
{\bf white dots} represent all the marked points on the curve~$C$
different from the critical points and from $x_1, \dots, x_{n+p}$.
We will not show them in the figures when it is not necessary.

\subsection{Intersection numbers in genus~0}
\label{Ssec:genus0}

In our computations we will need the value of the bracket
$$
\left< \tau_{1,a_1} \tau_{0,a_2} \dots \tau_{0,a_n} \right> = 
\int_{\oM^{1/r}_{g;a_1, \dots, a_n}} c_W \psi_1
$$
for $\sum a_i = (n-1)r$.
The {\em topological recursion relation} 
expresses $\psi_1$ as a sum of divisors:

\setlength{\unitlength}{1em}
\begin{equation} \label{Eq:toprec}
\psi_1 =  \epsfbox[0 15 155 30]{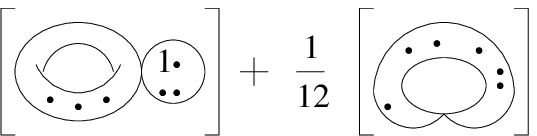} \; .
\end{equation}
\hspace{1em}

We must integrate Witten's class over each of these divisors.
The first divisor contributes $0$~\footnote{Indeed, the integral
of $c_W$ over this divisor includes as a factor 
the integral of $c_W$ over a genus~1 moduli space, which vanishes 
for reasons of dimension.} \label{snoska}
while the second one contributes
$$
\frac1{24} \sum_{a'+a'' = r-2} 
\left< 
\tau_{0,a_1} \dots \tau_{0,a_n} \tau_{0,a'} \tau_{0,a''}
\right>_0.
$$
This is a combination of integrals of Witten's class (without
$\psi$-classes) over genus~0 moduli spaces, and we are now going
to determine its value. For simplicity, in this section we
will use nonstandard notation for the bracket:
$\left<a_1, \dots, a_n \right>$ instead of
$\left<\tau_{0,a_1} \dots \tau_{0,a_n} \right>$.

\begin{proposition} \label{Prop:Loop}
For any $m \leq r-2$ and for any $x_1, \dots, x_n$,
$0 \leq x_i \leq r-2$, $\sum x_i = nr-m-2$, we have
$$
\sum_{a+b=m} \left<a,b, x_1, \dots, x_n \right> =
\frac{(n-1)!}{r^{n-1}} \prod_{i=1}^n (r-1-x_i).
$$
\end{proposition}

\paragraph{Proof.} The proof uses the known initial values of
the bracket:
$$
\begin{array}{rclcl}
\left<a_1,a_2,a_3\right> &=& 1& \mbox{for} &\sum a_i=r-2,\\
\left<a_1,a_2,a_3,a_4\right> &=& \frac1r \min(a_i, r-1-a_i)&
 \mbox{for}& \sum a_i =2r-2,
\end{array}
$$
and the WDVV equation. First note that the proposition is
true for $n=1$ (it says that $\sum_{a+b=m} 1 = m+1$).

Now we assume $n\geq 2$ and proceed by induction on~$m$.

If $m=0$ then the only possible bracket is $\left<0,0,r-2\right>$.

Suppose that the proposition is true up to some $m$. Let us apply the
WDVV equation to the correlators containing $1, a, b, x_1, \dots,
x_n$, the four distinguished points being $1$, $a$, $b$, and~$x_1$.
In other words, we consider all possible degenerations of the sphere into two
components such that $1$ and $a$ lie on one component, while
$b$ and $x_1$ lie on the other one, and then we swap the four points.
We obtain the following equality. (The summation over
$a+b=m$ is implicitly assumed; a hat means that the symbol
is skipped; underlined symbols are not in the original list,
they appear at the node of the degenerate sphere in the WDVV
formula.)
$$
\left<1,a,\underline{r-3-a}\right>
\left<\underline{a+1},b,x_1, \dots, x_n\right> \hspace{6.67cm} (A)
$$
$$
+ \sum_{i \not= 1} \left<1,a,x_i,\underline{2r-3-a-x_i}\right>
\left<\underline{a+1+x_i-r},b,x_1, \dots, \widehat{x_i}, \dots, x_n\right> 
\hspace{1cm} (B)
$$
$$
= \left<1,x_1,\underline{r-3-x_1}\right>
\left<a,b,\underline{x_1+1},x_2, \dots, x_n\right> \hspace{5.15cm} (C)
$$
$$
+ \sum_{i \not= 1} \left<1,x_1,x_i,\underline{2r-3-x_1-x_i}\right>
\times \hspace{10cm}
$$
$$
\hspace{4.5cm}\times \left<a,b,\underline{x_1+x_i+1-r},x_2, 
\dots, \widehat{x_i}, \dots, x_n\right>.
\hspace{.95cm} (D)
$$

\medskip

The term $(A)$ is the sum over $(a+1)+b = m+1$ that we want to
determine. The missing case $a+1=0$ does not matter, because
for $n \geq 2$ a single zero entry makes a bracket vanish.

The term $(B)$ can be evaluated by the induction assumption.
It vanishes if $a+1+x_i < r$, while for $a+1+x_i \geq r$
we obtain a sum over $(a+1+x_i)+b = m+x_i-r < m$. Thus we have
$$
(B) = \frac1r \cdot \frac{(n-2)!}{r^{n-2}} 
\prod_{j \not= i} (r-1-x_j).
$$

The term $(C)$ can, once again, be evaluated by the induction 
assumption:
$$
(C) = \frac{(n-1)!}{r^{n-1}} 
(r-2-x_1) \prod_{j \not= 1} (r-1-x_j).
$$

The last term $(D)$ vanishes if $x_1+x_i+1 <r$. Luckily, it turns out
that for any $i >1$, we have $x_1+x_i+1 \geq r$. Indeed, if
$a+b=m \leq r-3$ and $x_1+x_j \leq r-2$, then the sum of the
$n-2$ remaining terms equals $nr-m-2-x_1-x_i \geq (n-2)r+3$,
which is impossible since each of them equals at most $r-2$.
Therefore we have
$$
(D) = \frac1r \cdot \frac{(n-2)!}{r^{n-2}} 
[(r-1-x_1)+(r-1-x_i)]
\prod_{j \not= 1,i} (r-1-x_j).
$$
We deduce that
$$
(A) = (C)+(D)-(B) = \frac{(n-1)!}{r^{n-1}} \prod_j (r-1-x_j).
$$
\qed

\begin{remark}
By looking through the proof carefully one can check that the 
formula given in the proposition actually holds for
$m \leq r$, except if $n=1$.
\end{remark}

\section{Computations with double ramification cycles}
\label{Sec:formula}

In this section we prove Theorem~\ref{Thm:formula}. By
algebro-geometric arguments we find several relations
binding the values of the brackets involved in the theorem. In 
Section~\ref{Ssec:Relations} we write down these relations
and prove that they suffice to determine the values of the
bracket in all cases. We prove the relations in 
Section~\ref{Ssec:ProofRel}.

\subsection{The relations}
\label{Ssec:Relations}

Denote by 
$$
B = \left< \tau_{1, a_1} \tau_{0,a_2} \dots, \tau_{0, a_n} \right>_{g=1}
= \frac1{24} \, \frac{(n-1)!}{r^{n-1}}
\prod_{i=1}^n (r-1-a_i)
$$
(see Section~\ref{Ssec:genus0}).

Then the following relations hold.

\paragraph{Relation 1.}
$$
(k_1+1)(n_+ + n_- + 1) B \; = 
$$
$$
-(k_1 + n_+ + n_- + 1) 
\left<
\begin{array}{ccc|c|ccc}
k_1 & \dots & k_{n_+} & 0 \; \dots \;  0 & \tk_1 & \dots & \tk_{n_-} \\
a_1 & \dots & \dots & \dots & \dots & \dots & a_n
\end{array}
\right>
$$
$$
- \; \sum_{i=2}^{n_+} (k_i-1)
\left<
\begin{array}{cccccc|c|ccc}
k_1+1 & k_2 & \dots & k_i-1 & \dots & k_{n_+} 
& 0 \; \dots \; 0 & \tk_1 & \dots & \tk_{n_-} \\
a_1 &&& \dots && \dots & \dots & \dots & \dots & a_n
\end{array}
\right>
$$
$$
+ \; \sum_{i=1}^{n_-} (\tk_i+1)
\left<
\begin{array}{cccc|c|ccccc}
k_1+1 & k_2 & \dots & k_{n_+} 
& 0 \;  \dots \;  0 & \tk_1 & \dots & \tk_i+1 & \dots & \tk_{n_-} \\
a_1 & \dots && \dots & \dots & \dots && \dots && a_n
\end{array}
\right> .
$$

\paragraph{Relation 2.}
$$
(k_1+1) B \; = 
- \; \left<
\begin{array}{ccc|c|ccc}
k_1 & \dots & k_{n_+} & 
\begin{picture}(3.1,1)
\put(0,0){$\overbrace{0 \; \dots \;  0}^{n_0}$}
\end{picture}
& \tk_1 & \dots & \tk_{n_-} \\
a_1 & \dots & \dots & \dots &  \dots & \dots & a_n
\end{array}
\right>
$$
\vspace{1em}
$$
+ \; \left<
\begin{array}{cccc|c|cccc}
k_1 +1 & k_2 &\dots & k_{n_+} & 
\begin{picture}(3.1,1)
\put(0,0){$\overbrace{0 \; \dots \;  0}^{n_0-1}$}
\end{picture}
& 1 & \tk_1 & \dots & \tk_{n_-} \\
a_1 & \dots & \dots & \dots &  \dots & \dots & \dots & \dots & a_n
\end{array}
\right> \; .
$$

\paragraph{Relation 3.}
$$
\left<
\begin{array}{c|c|c}
1 & 0 \; \dots \; 0 & 1\\
a_1 & \dots & a_n
\end{array}
\right> = 0.
$$

\bigskip

To these relations we may add the simple observation
that the brackets are invariant under renumberings of
the marked points and under a simultaneous change of
sign of all $k_i$'s.

\bigskip

One can check by direct computations that these relations are
compatible with the expression of the bracket given in 
Theorem~\ref{Thm:formula}. 

\begin{lemma}
Relations~1, 2, and~3 determine the values of all 
brackets unambiguously.
\end{lemma}

\paragraph{Proof.} The lemma is proved by induction on
the number of nonzero entries $n_+ + n_-$ in a bracket,
the base of induction being given by Relation~3.

Consider a bracket whose value we would
like to determine.

First case: suppose that one of the $\tk_i$'s is equal to~1 and
one of the $k_i$'s is greater than~1 (or vice versa). 
Then we take our bracket
to be the last term in Relation~2 and replace it by the
first term (plus a know multiple of~$B$). We have decreased
the number of nonzero $k_i$'s in the bracket to be computed.

Second case: suppose that all the (nonzero) $k_i$'s and $\tk_i$'s equal~1.
In particular, $n_+ = n_-$. If $n_+=n_-=1$ the bracket vanishes
by Relation~3. If $n_+ = n_- \geq 2$, we take our bracket to be
the first term in Relation~1. Our bracket is than replaced
by a sum of brackets in which one of the $k_i$'s and one of the
$\tk_i$'s are equal to~2, while all the others are equal to~1.
Such brackets fall into the first case.

Third case: suppose that all of the $k_i$'s and $\tk_i$'s are greater
than~1. Then we take our bracket to be the first term of the
second sum in Relation~1 and replace it by the sum of the
other terms. As a result we obtain a combination
of brackets in all of which the number $\tk_1$ is smaller than
in the initial bracket. Repeating the same operation for each
bracket we are sure that $\tk_1$ will decrease by~1 with every
step. Thus after a finite number of steps we will end up with
a collection of brackets all of which fall into the first
two cases. \qed

\bigskip

Thus the relations determine the values of the brackets
unambiguously. Since they are compatible with the
expression of Theorem~\ref{Thm:formula}, the theorem
will be completely proved once we will have proved the
relations.

\subsection{Proof of the relations}
\label{Ssec:ProofRel}

First of all, note that Relation~3 is obvious. The
bracket in this relation denotes an integral over
an empty space, therefore it vanishes. Thus we
only need to prove Relations~1 and~2.

\begin{figure}[h]
\begin{center}
\setlength{\unitlength}{1em}
\begin{picture}(11,9)(0,-4)
\put(0,0){$\oM_{0,n+3}$}
\put(.5,6){$\oA$}
\put(1,5.3){\vector(0,-1){3.5}}
\put(6.5,0){$\oM_{1,n}$}
\put(6.5,6){$\oM_{1,(n+1)(K+1)}$}
\put(7.2,5.3){\vector(0,-1){3.5}}
\put(2.5,6.4){\vector(1,0){3}}
\put(2.5,5.3){\vector(1,-1){3.5}}
\put(0,3.5){$f$}
\put(4.5,3.5){$h$}
\put(4,6.8){$j$}
\put(7.7,3.5){$\pi$}
\end{picture}
\hspace{1em}
\epsfbox{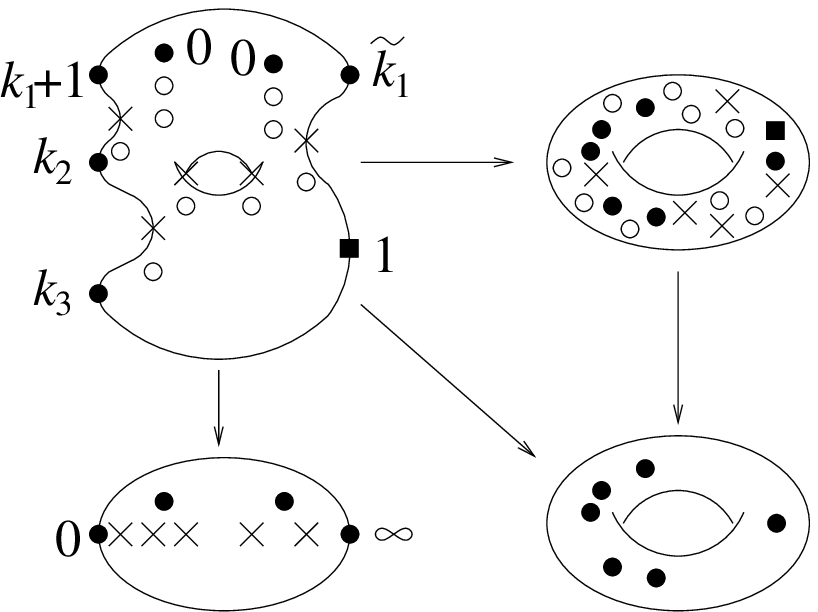}
\caption{Four moduli spaces.}
\label{Fig:diagramma}
\end{center}
\end{figure}

Consider the commutative
diagram involving four moduli spaces
shown in Figure~\ref{Fig:diagramma}.

Here $\oA = \oA(k_1+1, \dots, k_n,-1)$ is 
the DR-space of genus~1 admissible coverings with ramification
types $(k_1+1, k_2, \dots, k_{n_+})$ and 
$(1, \tk_1, \dots, \tk_{n_-})$ over $0$ and $\infty$. 
There are also $n_0$ additional marked points
on the source curve.

The projection $\pi$ forgets all marked points except the
round black dots and stabilizes the curve.

In this example, the DR-cycle coincides with the whole
moduli space $\oM_{1,n}$.

The main technique of this section is to compute the integral
$$
I = \int\limits_{\oA} h^*(c_W) \psi_1
$$
by three different methods and to use the relations thus 
obtained. The three methods can be summarized as follows.

\noindent
{\bf A.} Express the $\psi$-class on $\oM_{0,n+3}$
as a sum of boundary divisors. By Ionel's lemma
(Lemma~\ref{Lem:Ionel}),
the preimages under $f$ of these boundary divisors
represent the class~$\psi_1$ on~$\oA$. The images of
these divisors under~$h$ turn out to
to be DR-cycles. Therefore the integrals
of Witten's class over these divisors are some values
of the bracket of Definition~\ref{Def:divisor}.

\noindent
{\bf B.} Same as (A) with another expression of the 
$\psi$-class in terms of boundary divisors.

\noindent
{\bf C.} In Section~\ref{Ssec:genus0} we computed the integral
$$
B = \int_{\oM_{1,n}} c_W \psi_1
$$
analogous to $I$. Now, the classes $\psi_1$ on $\oM_{1,n}$
and on $\oA$ differ by a sum of boundary divisors. 
The images of these divisors under~$h$ turn out, once
again, to be DR-cycles, and integrating
Witten's class upon them we obtain a linear combination of brackets.

\paragraph{A.} On the space $\oM_{0,n+3}$ in
Figure~\ref{Fig:diagramma}, consider the class $\psi$
at the point labeled by $0$. This class is equal to
a sum of boundary divisors that we will represent as
\setlength{\unitlength}{1em}
\begin{equation} \label{Eq:psi1}
\psi = \epsfbox[0 15 75 50]{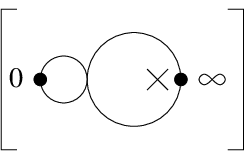} .
\end{equation}

\vspace{1em}

\noindent
The picture represents a sum of
boundary divisors where the sphere splits into two components.
The point labeled by $0$ is on the first component,
while the point labeled by $\infty$ together with
some {\em chosen} critical value (say, the first one)
is on the second component. The other marked points
on the sphere can be distributed arbitrarily between
the components.

Consider the preimages under $f$ of these boundary strata.

\begin{figure}[h]
\begin{center}
\
\begin{tabular}{lll}
\epsfbox{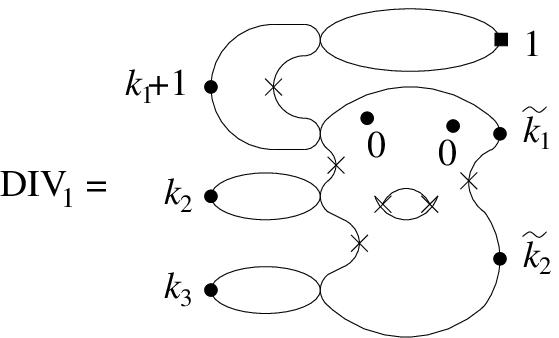} 
& \hspace{1em} &
\epsfbox{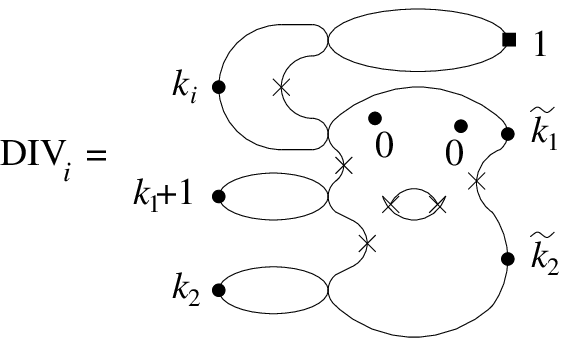} \\
\\
\epsfbox{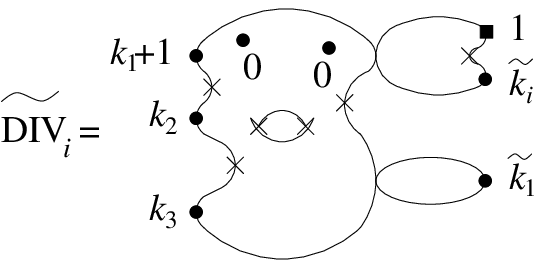} && \epsfbox{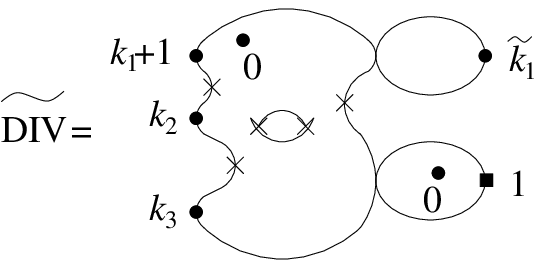}
\end{tabular}

\caption{Divisors in $\oA$.}
\label{Fig:DIV}
\end{center}
\end{figure}

\begin{lemma} \label{Lem:Adiv}
Among the preimages in $\oA$ of the divisors~{\rm (\ref{Eq:psi1})}
consider those on which the integral of $h^*(c_W)$ does not vanish.
These divisors are $\DIV_1$, $\DIV_i$, and $\tDIV_i$ in 
Figure~\ref{Fig:DIV}.
\end{lemma}

It should be understood that each picture actually represents
a {\em generic} admissible covering lying in the divisor.
In other words, the divisor is the closure of the set of
admissible coverings with the topological structure shown
in the figure.

\paragraph{Proof of the lemma.} We consider an irreducible
component of the preimage of~(\ref{Eq:psi1}) and reason in
terms of the generic admissible covering $\varphi: C \rightarrow S$
in this component.

First suppose that the curve $C$ does not have a toric component.
Then it looks as a ring of spheres with several ``tails'':
\begin{center}
\
\epsfbox{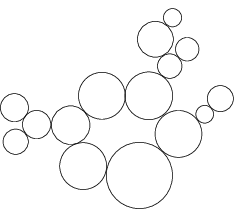}
\end{center}
Each sphere of the ring contains at least~2 nodes. Thus the
restriction of $\varphi$ to any such sphere is of degree
at least~2. Which implies that there is at least one
round black dot on each sphere of the ring (there is
at least one preimage of $0$ or $\infty$ on each component,
and the square black dot alone, being a simple preimage, is
not enough for $\varphi$ to have degree~2). It follows that
none of the spheres of the ring is contracted under the map~$h$
to $\oM_{1,n}$. Now, the number of spheres in the ring is even (those
that contain a preimage of~0 alternate with those that
contain a preimage of~$\infty$). Thus the curve $C$ retains
at least~2 nodes after the projection to $\oM_{1,n}$.
Which means that the codimension of the image of the corresponding divisor
under~$h$ is at least~2, so the integral of Witten's class on it vanishes.
The conclusion is that we only need to consider curves~$C$
with a toric component.

Suppose $C$ has a toric component with several ``tails'':
\begin{center}
\
\epsfbox{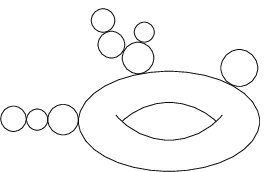}
\end{center}
If one of the tails contains more than one round black dot,
than the projection $h$ of such a curve to $\oM_{1,n}$
looks like:
\begin{center}
\
\epsfbox{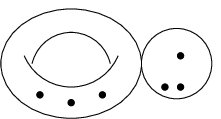}
\end{center}
The integral of Witten's class over the divisor of such curves
vanishes. The conclusion is that each ``tail'' contains at most
one round black dot.

Note that every component of $C$ contains at least one black dot
(round or square): either a preimage of~0 or of~$\infty$.
It follows that each
tail (except perhaps one) is composed of exactly~1 sphere containing
exactly~1 black dot. The exceptional tail can be composed of~1 or~2
spheres containing a round black dot and a square black dot. 
On each simple tail the function $\varphi$ has the form $z \mapsto z^k$
for some positive integer~$k$. 

Actually, an exceptional tail is bound to exist. Indeed,
suppose there is no exceptional tail, i.e., every tail
is composed of~1 sphere with a unique black dot. The image
curve~$S$ has~2 components, and in the case we consider,
the preimage of one of them is the toric component of $C$,
while the preimage of the other is the union of the tails.
Then the latter component of~$S$ contains no marked points
except $0$ or $\infty$. This is impossible, because~$S$
is stable.

If the exceptional tail is composed of~1 sphere, then this sphere
contains the black square dot and some dot $\tk_i$. This gives
the divisor $\tDIV_i$. If the exceptional tail is composed of~2
spheres, then one of them contains the square black dot, and the
other one some dot $k_i$ or $k_1+1$. The sphere with the
square dot can contain at most one node because the degree
of~$\varphi$ on it equals~1. Thus we obtain the divisors
$\DIV_1$ and $\DIV_i$.

The lemma is proved. \qed

Now, the integral $I$ is the sum of integrals of Witten's
class over the divisors $\DIV_1$, $\DIV_i$, and $\tDIV_i$, which we
will now determine.

\begin{lemma} \label{Lem:Acoefs}
The contributions to $I$ of the divisors $\DIV_1$, $\DIV_i$, and
$\tDIV_i$ equal
$$
\begin{array}{rl}
\frac{\displaystyle (n_++n_-)! (K-1)!^{n_++n_-+1} K!^{n_0}}
{\displaystyle k_1+1}
\;\; \cdot & 
k_1 (n_++n_-) \;\; \times\\
\\
\multicolumn{2}{r}{
\times
\left<
\begin{array}{ccc|c|ccc}
k_1 & \dots & k_{n_+} & 0 \; \dots \;  0 & \tk_1 & \dots & \tk_{n_-} \\
a_1 & \dots & \dots & \dots & \dots & \dots & a_n
\end{array}
\right>,
}\\
\\
\frac{\displaystyle (n_++n_-)! (K-1)!^{n_++n_-+1} K!^{n_0}}
{\displaystyle k_1+1}
\;\; \cdot &
(k_i-1) (n_++n_-) \;\;
\times\\
\\
\multicolumn{2}{r}{
\times
\left<
\begin{array}{cccccc|c|ccc}
k_1+1 & k_2 & \dots & k_i-1 & \dots & k_{n_+} 
& 0 \; \dots \; 0 & \tk_1 & \dots & \tk_{n_-} \\
a_1 &&& \dots && \dots & \dots & \dots & \dots & a_n
\end{array}
\right>,
}\\
\\
\frac{\displaystyle (n_++n_-)! (K-1)!^{n_++n_-+1} K!^{n_0}}
{\displaystyle k_1+1}
\;\; \cdot & 
(\tk_i+1) \;\; \times\\
\\
\multicolumn{2}{r}{\qquad \qquad
\times
\left<
\begin{array}{cccc|c|ccccc}
k_1+1 & k_2 & \dots & k_{n_+} 
& 0 \;  \dots \;  0 & \tk_1 & \dots & \tk_i+1 & \dots & \tk_{n_-} \\
a_1 & \dots && \dots & \dots & \dots && \dots && a_n
\end{array}
\right>,}
\end{array}
$$
respectively, where $K = \sum k_i = \sum \tk_i$.
\end{lemma}

\paragraph{Proof.}
The fundamental classes of the divisors $\DIV_1$, $\DIV_i$, and
$\tDIV_i$ project to multiples of double ramification
divisors in $\oM_{1,n}$ under the map~$h$. Thus we can
integrate Witten's class over these double ramification
divisors, which explains the brackets that appear in the
answers.

The coefficients in front of these brackets 
arise as products of three factors:
(i)~the transversal multiplicity of the map $f$ on the
divisor (because we need pull-backs under~$f$ of homology classes
rather than geometric pull-backs);
(ii)~the degree of $h$ on the divisor;
(iii)~the factor $1/(k_1+1)$ coming from Ionel's lemma.
For example, for $\DIV_i$, these factors are:

$$
{\rm (i)} \quad
(k_1+1) k_2 \dots k_{i-1} (k_i -1) k_{i+1} \dots k_{n_+}.
$$
These factors come from the number of ways to resolve the
nodes of the source curve.
$$
{\rm (ii)} \quad
\frac{(n_++n_-)(n_++n_-)! (K-1)!^{n_++n_-+1} K!^{n_0}}
{(k_1+1) k_2 \dots k_{i-1} k_{i+1} \dots k_{n_+}}.
$$
The factor $(n_++n_-)(n_++n_-)!$ is the number of ways to
number the ramification points on the two-component
genus~0 image curve, taking into account that the
first ramification point lies on the component of $\infty$.
$(K-1)!^{n_++n_-+1}$ is the number of ways to number the
white (noncritical) preimages of the ramification points.
$K!^{n_0}$ is the number of ways to number the white preimages
of the black marked points different from $0$ and~$\infty$.
The denominator comes from the spheres on the left
of the picture of $\DIV_i$: the restriction of the
admissible covering map $\varphi$ on these spheres
has the form $z^k$, and thus $z \mapsto (-1)^{1/k}z$
gives a renumbering of the white dots on such a sphere
equivalent to the initial numbering.

Multiplying these factors (without forgetting $1/(k_1+1)$)
we obtain the coefficient for $\DIV_i$ as claimed in the
lemma. The computations for the other divisors are analogous.
\qed

\paragraph{B.}
The second way of computing the integral~$I$ is not very different
from the first one. This time we start with a different 
presentation of the $\psi$-class at the point labeled by~$0$
on $\oM_{0,n+3}$:
\setlength{\unitlength}{1em}
\begin{equation} \label{Eq:psi2}
\psi = \epsfbox[0 15 75 50]{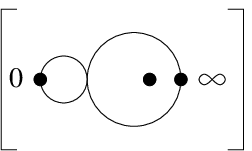} .
\end{equation}

\vspace{1em}

\noindent
The picture represents a sum of
boundary divisors where the sphere splits into two components.
The point labeled by~$0$ is on the first component,
while the point labeled by~$\infty$ together with
another {\em chosen} black marked point (say, the first one)
is on the second component.

\begin{lemma} \label{Lem:Bdiv}
Among the preimages in $\oA$ of the divisors~{\rm (\ref{Eq:psi2})}
consider those on which the integral of $c_W$ does not vanish.
These divisors are $\DIV_1$, $\DIV_i$, and $\tDIV$ 
in Figure~\ref{Fig:DIV}.
\end{lemma}

\paragraph{Proof.} First suppose there are no points
marked with crosses on the component of~$\infty$ on the
image curve. This means that the preimages of this component
are spheres containing only black marked points. However, 
all these components must be contracted under the map~$h$. (Otherwise the
integral of $c_W$ over this divisor vanishes.) Therefore, actually,
every sphere contains exactly one black point, except one
sphere that contains a round black dot and a square black dot.
This is the divisor $\tDIV$.

If the component of $\infty$ in the image curve contains
at least one cross-marked point, then we have reduced the
problem to the situation of Lemma~\ref{Lem:Adiv}. Thus the
only divisors that can give a nonzero contribution
are $\DIV_1$, $\DIV_i$, and $\tDIV_i$. However, actually
the divisors $\tDIV_i$ do not appear as preimages of
the divisors~(\ref{Eq:psi2}) because in this case the
component of $\infty$ in the image contains no black dots
other than $\infty$. \qed

\begin{lemma} \label{Lem:Bcoefs}
The contributions to $I$ of the divisors $\DIV_1$, $\DIV_i$,
and $\tDIV$ equal
$$
\begin{array}{rl}
\frac{\displaystyle (n_++n_-)! (K-1)!^{n_++n_-+1} K!^{n_0}}
{\displaystyle k_1+1}
\;\; \cdot & 
k_1 (n_++n_-+1) \;\; \times\\
\\
\multicolumn{2}{r}{
\times
\left<
\begin{array}{ccc|c|ccc}
k_1 & \dots & k_{n_+} & 0 \; \dots \;  0 & \tk_1 & \dots & \tk_{n_-} \\
a_1 & \dots & \dots & \dots & \dots & \dots & a_n
\end{array}
\right>,
}\\
\\
\frac{\displaystyle (n_++n_-)! (K-1)!^{n_++n_-+1} K!^{n_0}}
{\displaystyle k_1+1}
\;\; \cdot &
(k_i-1) (n_++n_-+1) \;\;
\times\\
\\
\multicolumn{2}{r}{
\times
\left<
\begin{array}{cccccc|c|ccc}
k_1+1 & k_2 & \dots & k_i-1 & \dots & k_{n_+} 
& 0 \; \dots \; 0 & \tk_1 & \dots & \tk_{n_-} \\
a_1 &&& \dots && \dots & \dots & \dots & \dots & a_n
\end{array}
\right>,
}\\
\\
\frac{\displaystyle (n_++n_-)! (K-1)!^{n_++n_-+1} K!^{n_0}}
{\displaystyle k_1+1}
\;\; \cdot & 
(n_++n_-+1) \;\; \times\\
\rule{0em}{1.3em}\\
\multicolumn{2}{r}{\qquad \qquad
\times
\left<
\begin{array}{cccc|c|cccc}
k_1 +1 & k_2 &\dots & k_{n_+} & 
\begin{picture}(3.1,1)
\put(0,0){$\overbrace{0 \; \dots \;  0}^{n_0-1}$}
\end{picture}
& 1 & \tk_1 & \dots & \tk_{n_-} \\
a_1 & \dots & \dots & \dots &  \dots & \dots & \dots & \dots & a_n
\end{array}
\right>,}
\end{array}
$$
respectively, where $K = \sum k_i = \sum \tk_i$.
\end{lemma}

The proof is analogous to that of Lemma~\ref{Lem:Acoefs}. \qed

\paragraph{C.}
Now we are going to evaluate
$$
I = \int_{\oA} h^*(c_W) \psi_1(\oA)
$$
using the value
$$
B = \int_{\oM_{1,n}} c_W \psi_1(\oM_{1,n}).
$$
If the class $\psi_1$ on $\oA$ were a pull-back of the
class $\psi_1$ on $\oM_{1,n}$, then we would simply have
$I = B \cdot \deg h$. However this equality is actually
incorrect, because there is difference between
$h^*(\psi_1(\oM_{1,n}))$ and $\psi_1(\oA)$.

\begin{lemma} \label{Lem:Cdiv}
The difference $\psi_1(\oA) - h^*(\psi_1(\oM_{1,n})$ can be represented
as the sum of divisors
$$
k_2 \dots k_{n_+} \; \DIV_1 \; + \; 
\sum_{i=2}^{n_+} k_2 \dots k_{i-1} (k_i -1) k_{i+1} \dots k_{n_+} \; \DIV_i
$$
{\rm (}the divisors being shown in
Figure~{\rm \ref{Fig:DIV})} and, in addition, 
some other divisors on which the integral of $c_W$ vanishes.
\end{lemma}

\paragraph{Proof.} There is a subtlety in finding the difference 
$\psi_1(\oA) - h^*(\psi_1(\oM_{1,n})$. First consider the
difference $\psi_1(\oM_{1,(n+1)(K+1)}) - \pi^*(\psi_1(\oM_{1,n})$
in the upper right moduli space in Figure~\ref{Fig:diagramma}.
It is given by the divisor $D_1$ of all stable curves on which the first
black marked point ($k_1+1$) is situated on a component contracted
by~$\pi$. Now we consider the intersection of this divisor
with the image $j(\oA)$. It turns out that this intersection is
not necessarily transversal, and its multiplicity gives the 
coefficients that appear in the formulation of the lemma.

Now we take the pull-back of the intersection by~$j$.
Let us consider an irreducible component of this pull-back
$$
j^* \Bigl(
\psi_1(\oM_{1,(n+1)(K+1)}) - \pi^*(\psi_1(\oM_{1,n}))
\Bigr)
$$
and a generic admissible covering $\varphi$ in this component.

The image curve of $\varphi$ in $\oM_{0,n+3}$ has a unique
node (otherwise the codimension of such a component in $\oA$
would be at least~2). 

Moreover, this node separates~$0$
and~$\infty$. Indeed, otherwise the component containing
the point $k_1+1$ also contains some preimage of $\infty$.
Since this component is contracted by $h$, the only
possible preimage is the square black dot (and there
are no other preimages of $\infty$). Thus the degree of $\varphi$
on this component equals~1. But this implies $k_1+1=1$, which
is impossible.

On the component of $\infty$ there exists at least one cross
marked point. Indeed, otherwise there is no ramification over this
component, so all the ramification is over the component containing~0.
Then the point $k_1+1$ lies on the torical component, which is
impossible, since this component should be contracted by~$h$.

Thus we have reduced the problem to the situation of
Lemma~\ref{Lem:Adiv} with the additional restriction that
the component containing the point $k_1+1$ should be
contracted by~$h$. This restriction excludes the divisors $\tDIV_i$
but allows the divisors $\DIV_1$ and $\DIV_i$. 

It remains to find the multiplicity of the intersection 
$j(\oA) \cap D_1$ along $j(\DIV_1)$ and $j(\DIV_i)$. 
This multiplicity is easily seen to be equal to the
product of indices of the nodes that are desingularized as a 
generic point of the intersection moves to a generic point of~$D_1$.
Such products of indices are precisely the coefficients stated in the
lemma. \qed

\begin{lemma} \label{Lem:Ccoefs}
The contributions to $I$ of the integral $B$ and that of the
divisors $\DIV_1$ and $\DIV_i$ equal
$$
\begin{array}{rl}
\frac{\displaystyle (n_++n_-)! (K-1)!^{n_++n_-+1} K!^{n_0}}
{\displaystyle k_1+1}
\;\; \cdot & 
(k_1+1) (n_++n_-+1) \;\; \cdot \;\; B, \qquad \qquad\\
\\
\frac{\displaystyle (n_++n_-)! (K-1)!^{n_++n_-+1} K!^{n_0}}
{\displaystyle k_1+1}
\;\; \cdot & 
k_1 (n_++n_-+1) \;\; \times\\
\\
\multicolumn{2}{r}{
\times
\left<
\begin{array}{ccc|c|ccc}
k_1 & \dots & k_{n_+} & 0 \; \dots \;  0 & \tk_1 & \dots & \tk_{n_-} \\
a_1 & \dots & \dots & \dots & \dots & \dots & a_n
\end{array}
\right>,
}\\
\\
\frac{\displaystyle (n_++n_-)! (K-1)!^{n_++n_-+1} K!^{n_0}}
{\displaystyle k_1+1}
\;\; \cdot &
(k_i-1) (n_++n_-+1) \;\;
\times\\
\\
\multicolumn{2}{r}{
\times
\left<
\begin{array}{cccccc|c|ccc}
k_1+1 & k_2 & \dots & k_i-1 & \dots & k_{n_+} 
& 0 \; \dots \; 0 & \tk_1 & \dots & \tk_{n_-} \\
a_1 &&& \dots && \dots & \dots & \dots & \dots & a_n
\end{array}
\right>,
}
\end{array}
$$
respectively, where $K = \sum k_i = \sum \tk_i$.
\end{lemma}

\paragraph{Proof.} The coefficient of $B$ is just the degree of~$h$.
The coefficients of the other two brackets are obtained as products
of two factors: 
(i)~the coefficients of the divisors appearing in Lemma~\ref{Lem:Cdiv},
(ii)~the degree of $h$ on the divisor.

For instance, for $\DIV_i$ these factors equal:
$$
{\rm (i)} \quad k_2 \dots k_{i-1} (k_i-1) k_{i+1} \dots k_{n_+}.
$$
$$
{\rm (ii)} \quad 
\frac{(n_++n_-+1)! (K-1)!^{n_++n_-+1} K!^{n_0}}
{(k_1+1) k_2 \dots k_{i-1} k_{i+1} \dots k_{n_+}}.
$$
The case of the divisor $\DIV_1$ is analogous. \qed

\bigskip

Thus we have established three expressions, A, B, and~C  
for the integral~$I$. Writing ${\rm A} - {\rm B} = 0$
and ${\rm A} - {\rm C} = 0$ we obtain Relations~1 and~2.

\section{An algorithm to compute Witten's intersection numbers}
\label{Sec:algorithm}

Here we present an algorithm for computing any number
$\left< \tau_{d_1, a_1} \dots \tau_{d_n, a_n} \right>$.
This algorithm is rather hard to implement on a computer
because it involves enumerating all possible degenerations
of an admissible covering satisfying some given properties. Thus
this section is best viewed as a constructive proof of
Theorem~\ref{Thm:algorithm}.

Suppose we are given $n$ nonnegative integers $d_1, \dots, d_n$.
Choose an integer $p$, $0 \leq p \leq g$ and 
$n$ numbers $k_1, \dots, k_n$ large enough (we will use the explicit
bound $|k_i| > \sum d_i$) such that $\sum k_i = p$. 
Consider the DR-space $\oA = \oA(k_1, \dots, k_n, -1, \dots, -1)$,
where the list ends with $p$ numbers~$-1$ so that the total sum is~0 as
it should be. 

As explained in Section~\ref{Ssec:Adm},
there is a projection $h:\oA \rightarrow \oM_{g,n}$ that 
sends $\oA$ to DR-cycle $\oD$ of codimension~$p$ in $\oM_{g,n}$.

We are going to compute the integral
\begin{equation} \label{Eq:integral}
\int_{\oA} h^*(c_W \, \psi_1^{d_1} \dots \psi_n^{d_n}) = 
\deg(h) \int_{\oD} c_W \, \psi_1^{d_1} \dots \psi_n^{d_n}.
\end{equation}
If $p=g$ this will give us the value of the bracket
$\left< \tau_{d_1,a_1} \dots \tau_{d_n,a_n} \right>$.
Note, however that our result is actually more general.

Our algorithm for computing Integral~(\ref{Eq:integral})
can be summed up as follows.

While the class $\psi_i$ on $\oM_{g,n}$ is, in general, not representable
by boundary divisors, its pull-back to $\oA$ can be easily represented
as a sum like that. At each step we will consider one
$\psi$-class and replace it by a
sum of divisors. Thus we will be reduced to computing the integral
of a product $c_W \cdot (\mbox{powers of } \psi\mbox{-classes})$ over some
divisors, the number of $\psi$-classes having decreased by~1.
Using the factorization rules for Witten's class, we will be
able to represent each integral over a boundary divisor as 
a product of analogous integrals over simpler DR-spaces.

\begin{proposition} \label{Prop:algorithm}
The integral~{\rm (\ref{Eq:integral})} can be effectively expressed
as a combination of analogous integrals with a smaller sum
$\sum d_i$.
\end{proposition}

\paragraph{Proof.} \

\noindent
{\bf 1.} Expressing the pull-back $h^*(\psi_i)$ as a
sum of boundary divisors of~$\oA$.

Consider the class $\psi_i = \psi_i(\oM_{g,n})$ 
on $\oM_{g,n}$ and the corresponding
class $\Psi_i = \psi_i(\oA)$ on $\oA$. The pull-back $h^*(\psi_i)$ to
$\oA$ can be represented as a sum of divisors in the following way.
First, the class $\Psi_i$
can, by Ionel's lemma (Lemma~\ref{Lem:Ionel}), be replaced
by the pull-back of the $\psi$-class at $0$ or at $\infty$
on $\oM_{0,n+p+2g}$. The latter class is easy to represent
by a sum of divisors, see Equation~(\ref{Eq:psi1}). 
Second, the difference $\Psi_i - h^*(\psi_i)$ is equal to the sum
(with certain coefficients) of divisors $D_i$ formed by the 
admissible coverings for which the component of the source curve
containing the $i$th marked point is contracted by~$h$
(cf.~Lemma~\ref{Lem:Cdiv}). 

\medskip

\noindent
{\bf 2.} The image curve of a generic admissible covering in every
boundary divisor has a unique node separating~$0$ and~$\infty$.

Let us consider an irreducible component of one of the above divisors
and a generic admissible covering $\varphi$ in this component.
The image curve of $\varphi$ in $\oM_{0,n+p+2g}$ has a unique
node because otherwise the codimension of such a component in $\oA$
would be at least~2. 

Let us prove that this node separates~$0$
and~$\infty$. Indeed, this is obvious by construction 
for the divisors involved in the expression of $\Psi_i$.
As for the divisors representing the difference
$\Psi_i - h^*(\psi_i)$, suppose that the node on the image curve
does not separate $0$ and $\infty$. Consider the component 
of the source curve $C$ containing the $i$th black round marked
point. This compounent must contain both a preimage of $0$ and
a preimage of~$\infty$. On the other hand, it is, by construction,
contracted by $h$. Thus the only possibility is that it contains a
unique round black dot (a preimage of $0$) and one or several square
black dots (since these are forgotten by~$h$). But then the degree of
$\varphi$ on this component is equal to $|k_i|$ and, at the same time,
to the number of square black dots. Since we assumed that 
$|k_i| > \sum d_i$ this is impossible. 

\medskip

\noindent
{\bf 3.} Splitting the admissible covering.

Denote by $\DIV$ the boundary divisor of $\oA$ under
consideration.

\begin{figure}
\begin{center}
\
\epsfbox{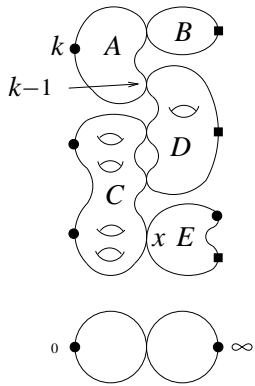}
\caption{A boundary divisor is isomorphic to a product of several
simpler DR-spaces.}
\label{Fig:split}
\end{center}
\end{figure}

Since a generic admissible covering in $\DIV$ can
have nontrivial ramifications only over $0$, $\infty$, and the
node of the image curve, we see that there are at most 2 
multiple ramification points on each component of the image curve. 
Thus, we can split the admissible covering according to the 
components of the source curve and 
obtain several simpler admissible coverings lying in
simpler DR-spaces. This is shown in Figure~\ref{Fig:split}
(the meaning of various markings will be explained later).
We will call the components of the source curve in a
generic admissible covering {\em parts}. In Figure~\ref{Fig:split}
the parts are denoted by $A$, $B$, $C$, $D$, and~$E$.

The smaller admissible coverings will 
have a lot of unnecessary white marked points. Consider,
for instance, a critical point (a cross) on the part~$D$
in Figure~\ref{Fig:split}. Its image is a ramification point 
(a cross) in the image curve. This cross has several simple
(white) preimages on the parts~$B$ and~$E$. Once we have
separated the parts we would like to forget these
useless white points. This is done in the following way:
(i)~forget those marked points on the image curve that
have only white preimages (together with these white
preimages), (ii)~contract the unstable components
of the image curve and their preimages in the source curve.
If the image curve happens to have only marked points, then
we simply ignore the corresponding part.
For example, in Figure~\ref{Fig:split} the part $B$
will be ignored.

Denote by $c$ the number of parts that
remain after that ($c=4$ in our example), and denote by
$\oA_1, \dots, \oA_c$ the corresponding DR-spaces. We
have constructed a map
$$
u: \DIV \rightarrow \oA_1 \times \dots \times \oA_c \, .
$$
This map sends the fundamental homology class of $\DIV$
to a (multiple of) the fundamental homology class of
the product in the image. Indeed, $u$ is an isomorphism
on the open part of $\DIV$ where ramification points
do not coincide.

\medskip

\noindent
{\bf 4.} Re-assigning the $\psi$-classes.

Recall that the map $h$ takes the source curve in
Figure~\ref{Fig:split}, forgets
all the marked points except the round black dots and stabilizes
the curve. We are interested in the preimages of the
$\psi$ classes under $h$. During the stabilization the marked points
from certain parts can land on some other parts.
In the figure, the $\psi$-class assigned to the point
labeled with $k$ (on the part~$A$) will land on the part~$D$.
Similarly, the $\psi$-class assigned to the round black dot
of the part~$E$ will land on the part~$C$. Indeed, for
{\em any} admissible covering lying in the boundary divisor
represented in the figure, the parts $A$, $B$, and $E$
are contracted by~$h$. Once we have re-assigned the $\psi$-classes
in that way, the preimages $h^*(\psi_i)$ on $\DIV$ coincide with
the classes $h^*(\psi_i)$ defined separately for each part.
Indeed, for instance, in our example, stabilizing the source
curve is equivalent to contracting the parts $A$, $B$, and $E$
and then stabilizing the parts~$C$ and~$D$.

\medskip

\noindent
{\bf 5.} Splitting the integral.

Now we can split the integral~(\ref{Eq:integral}) into a product
of similar integrals over $\oA_1$, \dots, $\oA_c$ with the 
$\psi$-classes assigned as explained in paragraph~4.

The total number of $\psi$-classes in these integrals
equals $\sum d_i -1$, because we replaced one of the 
$\psi$-classes by the boundary divisors.

This completes the proof of the proposition.
\qed

\bigskip

Proposition~\ref{Prop:algorithm} constitutes the recursive
step of our algorithm. To make things precise we must add two 
comments.

{\bf a)} The condition $|k_i| > \sum d_i$ is easily seen
to be still satisfied for the smaller DR-spaces. Indeed,
An index $k$ of a zero or a pole can decrease only
by ``annihilating'' one or several black squares. For
instance, in Figure~\ref{Fig:split} the initial index $k$
has become equal to $k-1$ on the component $D$ by annihilating
one square black dots. 

{\bf b)} The restriction of Witten's class (more precisely,
of $h^*(c_W)$) on the components
is obtained by using the factorization rules. This involves
some choices or remainders modulo~$r$. For example, in 
Figure~\ref{Fig:split}, we will have to sum over $r-1$
possibilities of the remainders at the nodes connecting
$C$ to~$D$. The remainders assigned to the sqaure dots are equal
to~0 since they are forgotten under the map~$h$.

\paragraph{Proof of Theorem~\ref{Thm:algorithm}.}
Assume we want to find the value of an integral over $\oM_{g,n}$
involving Witten's class and powers of the $\psi$-classes.
First, using the above lemmas we can get rid of the
$\psi$-classes and reduce the problem to computing the
integral of Witten's class over DR-spaces. Comparing the
degree of Witten's class 
$$
\deg c_W = \frac{(r-2)(g-1)+\sum a_i}r \leq \frac{(r-2)(n+g-1)}r
$$
to the dimension of a DR-cycle
$$
\dim = 3g-3+n-p \quad (p \leq g),
$$
we see that such an integral can be nonzero only
in 2 cases: either for genus~0 or for DR-cycles of codimension~1
in genus~1. The genus zero integrals are well-known, while
the case of double ramification divisors in genus~1 was treated
in Section~\ref{Sec:formula}. \qed

\begin{remark}
Our first goal was to compute the integral
of $c_W \psi_1^{d_1} \dots \psi_n^{d_n}$ over the moduli
spaces $\oM_{g,n}$, but we actually computed such integrals over all
DR-cycles satisfying the strange-looking condition
$|k_i| > \sum d_i$. This generalization is unavoidable
if we want to make the algorithm work. It is easy to see
that even if we start with an integral over a moduli space the 
integral can immediately lead us to integrals over DR-cycles.
\end{remark}

\end{document}